# In Orbifolds, Small Isoperimetric Regions are Small Balls

Frank Morgan
Department of Mathematics and Statistics
Williams College
Williamstown, MA 01267
Frank.Morgan@williams.edu

## Abstract

In a compact orbifold, for small prescribed volume, an isoperimetric region is close to a small metric ball; in a Euclidean orbifold, it is a small metric ball.

## 1. Introduction

Even in smooth Riemannian manifolds M, there are relatively few examples of explicitly known regions which minimize perimeter for prescribed volume (see [M1, 13.2], [HHM]). One general result is that for M compact, for small volume, a perimeter-minimizing region is a nearly round ball where the scalar curvature is large ([K], [MJ, Thm. 2.2], or in 3D [Ros, Thm. 18], with [Dr]).

For singular ambients, one result in general dimensions is that in a smooth cone of positive Ricci curvature, which has just a single singularity at the apex, geodesic balls about the apex minimize perimeter [MR, Cor. 3.9]. Similarly in the surface of a convex polytope in **R**$^n$, with its stratified singular set, for small prescribed volume, geodesic balls about some vertex minimize perimeter [M2, Thm. 3.8]. Our Theorem 3.4 proves a similar result for orbifolds, which unlike the previous categories, are not topological manifolds. An orbifold is locally a Riemannian manifold modulo a finite group of isometries (see Definition 2.1). Corollary 3.3 concludes that in a Euclidean orbifold, for small volume, an isoperimetric region is a metric ball.

**The proof.** The proof of our main Theorem 3.4 has the following steps, sometimes following [M2] and predecessors, sometimes needing new arguments special to orbifolds. Let $R_\alpha$ be a sequence of isoperimetric regions in the orbifold $O$ with volumes approaching zero.

(1) After covering *O* by maps from unit lattice cubes in $\mathbf{R}^n$, Lemma 3.1 obtains a concentration of volume after [M1, 13.7].

(2) After covering *O* by geodesic balls modulo finite groups and isometrically embedding the balls in some Euclidean space, we obtain a limit of blow-ups of preimages $T_\alpha$ of the $R_\alpha$, nonzero by (1).

(3) After Almgren [Alm] and [M2, Lemma 3.5], we carefully show that $\partial T_\alpha$ locally minimizes area plus $C\Delta V$ among invariant surfaces.

(4) Using (3), generalizing "monotonicity" as in [M2, Lemma 3.6], we obtain a lower bound on the area of $\partial T_\alpha$ in any unit ball about one of its points and deduce that components of isoperimetric regions have small diameter.

(5) As in [MJ], it follows that for small prescribed volume, isoperimetric regions are close to being small metric balls.

**Acknowledgments.** I took up this work in response to a question from Antonio Ros, and I finished the proof of Theorem 3.4 while visiting Jaigyoung Choe at the Korean Institute for Advanced Study in November, 2006. I realized Propositions 3.5 and 3.6 at the February 2007 Levico workshop on "Questions about Geometric Measure Theory and Calculus of Variations," organized by R. Serapioni, F. Serra Cassano and I. Tamanini of Trento University and the Centro Internazionale per la Ricerca Matematica.

## 2. Orbifolds and Invariant Regions

**2.1. Definitions.** For $n \geq 1$, an nD (smooth Riemannian) *orbifold O* is a connected metric space locally isometric to a smooth nD Riemannian manifold M modulo a finite group G of isometries. Locally, the rectifiable currents of *O* (the generalized kD surfaces of geometric measure theory) are just the G-invariant rectifiable currents of M. (The locally rectifiable currents of M are defined by isometrically embedding M in some Euclidean space. The space and its topology are independent of the embedding. See [M1].)

The following lemma makes many basic geometric properties of orbifolds obvious.

**2.2. Structure Lemma.** *An nD orbifold is a smooth Riemannian manifold with boundary except for a smoothly stratified (n-2)D singular set.*

We mean that the singular set is an (n-2)D manifold except for a smoothly stratified (n-3)D singular set defined recursively down to isolated singular points.

*Proof by induction.* The case n = 1 is trivial. For n > 1, let *p* be a point in an nD orbifold *O*, which is locally isometric to a manifold M modulo a finite group of isometries G leaving a corresponding point *q* fixed. The tangent space modulo G is the cone over the unit sphere modulo

G, for which the lemma holds by induction. Therefore the lemma holds for the tangent space modulo G and hence for M/G, hence for *O*.

**2.3. Observation** [M2, Rmk. after 3.3]. Let G be a finite subgroup of $\mathbf{O}_n$. Since balls are uniquely perimeter minimizing in $\mathbf{R}^n$, G-invariant balls are uniquely minimizing among invariant regions. In other words, in $\mathbf{R}^n$ modulo a finite subgroup of $\mathbf{O}_n$, a perimeter-minimizing region for fixed volume is a metric ball about the image of a fixed point.

The following regularity theorem of [M3] applies indirectly to orbifolds, since surfaces in orbifolds locally correspond to invariant surfaces in manifolds.

**2.4. Theorem.** *In a smooth ($C^\infty$) Riemannian manifold M, among hypersurfaces oriented by unit normal invariant under a group G of isometries of M, perhaps with given boundary or homology or volume, suppose that S minimizes area. Then on the interior S is a smooth constant-mean-curvature hypersurface except for a singular set of codimension at least 7 in S.*

*Proof.* This theorem is stated and proved in [M3, Thm. 4.1], without full details for the case of a volume constraint, which we now provide. There is an $\Upsilon > 0$ such that a hypercone which is area-minimizing without volume constraint among invariant surfaces and which has density less than $\Upsilon$ is a hyperplane [M3, Lemma 4.4]. By Allard's regularity theorem [A, Sect. 8], S is regular ($C^{1,1/2}$) at all such points. By Federer [F, Lemma 2] the Hausdorff dimension of the singular set of S is less than or equal to the Hausdorff dimension of the singular set of some area-minimizing (tangent) hypercone in $\mathbf{R}^n$, which is at most n−8, i.e., codimension at least 7 in S. Higher smoothness follows by Schauder theory (see [M4, Prop. 3.5]).

## 3. Isoperimetric Regions in Orbifolds

Our main Theorem 3.4 follows from Proposition 3.2 and the methods of [MJ], which treats the easier case of manifolds. First we need a lemma on concentration of volume after [M1, 13.7]. For the lemma we will work in the category of metric spaces *O* which are smooth nD manifolds except for a set of nD Hausdorff measure 0. Regions and their boundaries are technically locally rectifiable currents in the regular part of *O* [M1].

**3.1. Lemma.** *For $n \geq 1$, let O be a metric space which is a smooth nD manifold except for a set of nD Hausdorff measure 0, covered by some finite number $k_1$ of distance-decreasing maps from disjoint open unit lattice cubes in $\mathbf{R}^n$ each of multiplicity bounded by $k_2$ which shrink (n−1)D or nD Hausdorff measure by at most some constant factor $C_1$. Then given C > 0 there is a constant $\delta > 0$ such that for any region R in O with topological boundary $\partial R$ and with nD and (n−1)D Hausdorff measure satisfying $|\partial R| \leq C|R|^{(n-1)/n}$, the image K of some relevant lattice subcube of volume between $2^{-n}|R|$ and $|R|$ satisfies*

(1) $$|R \cap K| \geq \delta |R|.$$

*Proof.* Subdivide the unit lattice cubes into congruent subcubes $J_i$ with volume between $|R|$ and $2^{-n}|R|$. Let $K_i$ denote the image of $J_i$ in $O$ and let $R_i$ denote the preimage of R in $J_i$. We may assume that $\delta < 1/2$ and $|R_i| \leq |J_i|/2$, since otherwise $|R \cap K_i| \geq |J_i|/2C_1 \geq |R|/2^{n+1}C_1$ and we are done. By the relative isoperimetric inequality [M1, 12.3(1)] for open cubes, there is a single constant $\gamma$ such that

$$|\partial R_i| \geq \gamma |R_i|^{(n-1)/n}.$$

Hence

$$|\partial R \cap K_i| \geq \gamma |R \cap K_i|^{(n-1)/n}/C_1 \geq \gamma |R \cap K_i| / C_1 \max|R \cap K_j|^{1/n}.$$

Since the total multiplicity of covering is at most $k = k_1 k_2$, summing over i yields

$$k|\partial R| \geq \gamma |R| / C_1 \max|R \cap K_j|^{1/n},$$

$$\max|R \cap K_j|^{1/n} \geq \gamma |R| / kC_1|\partial R| \geq \gamma |R|^{1/n} / kCC_1.$$

Hence some K satisfies
$$|R \cap K| \geq \delta |R|$$

with $\delta = (\gamma / kCC_1)^n$, as desired.

**3.2. Proposition.** *For $n \geq 2$, let O be a compact nD orbifold. For small prescribed volume, the sum of the diameters of the components of a perimeter-minimizing region is small.*

*Proof.* Fix finitely many Riemannian geodesic balls $M_i$ modulo finite groups $G_i$ of rotations about the center which cover $O$. We may assume that the $M_i$ are isometrically embedded in some fixed $\mathbf{R}^N$.

Let $R_\alpha$ be a sequence of perimeter-minimizing regions in $O$ with volume $|R_\alpha|$ small and approaching 0. Comparison with small geodesic balls about a regular point of $O$ shows that $|\partial R_\alpha| \leq C|R_\alpha|^{(n-1)/n}$. For each $R_\alpha$ use Lemma 3.1 to choose $K_\alpha$ with $2^{-n}|R_\alpha| \leq |K_\alpha| \leq |R_\alpha|$ and

$$|R_\alpha \cap K_\alpha| \geq \delta_0 |R_\alpha| \geq \delta_0 |K_\alpha|.$$

By taking a subsequence, we may assume that each $K_\alpha$ is contained well inside one fixed $M_i/G_i$, which we will now call M/G. Rescale $K_\alpha$ and M up to $1 = |R_\alpha| \geq |K_\alpha| \geq 2^{-n}$; denote the rescaling of M by $M_\alpha$, with extrinsic curvature going to zero. We may assume that $M_\alpha$ is tangent to $\mathbf{R}^n \times \{0\}$ at the origin at a point $p_\alpha$ in the preimage of $K_\alpha$.

By compactness [M1, 9.1], we may assume that the preimages $T_\alpha$ of $R_\alpha$ in $M_\alpha$ converge weakly to a region T with $|G| \geq |T| \geq 2^{-n}\delta_0 > 0$, lying in $\mathbf{R}^n \times \{0\}$.

Consider any g in G, and let $q_\alpha$ be its fixed point in $M_\alpha$ closest to the origin of $\mathbf{R}^N$. Unless the sequence $q_\alpha$ diverges to infinity, by taking a subsequence we may assume that $q_\alpha$ converges to a point q in $\mathbf{R}^n \times \{0\}$ and that the action of g on $M_\alpha$ converges to an isometry g of $\mathbf{R}^n \times \{0\}$ with fixed point q, in the sense that if points $x_\alpha$ in $M_\alpha$ converge to x in $\mathbf{R}^n\times\{0\}$, then $g(x_\alpha)$ converges to g(x), uniformly on compacts. Thus we obtain a limit action of some subgroup $G_1$ of G on $\mathbf{R}^n\times\{0\}$, which extends trivially to $\mathbf{R}^N$.

There are many points of $\partial T$ not on the cone of fixed points of $G_1$. Choose a small positive $\delta_1$ so that there are two $\delta_1$-balls $B_1$, $B_2$ in $\mathbf{R}^n\times\{0\}$ about such points of $\partial T$ with disjoint images under $G_1$ such that any third $\delta_1$-ball B is disjoint from the image under $G_1$ of one of them, say $B_1$. By the Gauss-Green-De Giorgi-Federer Theorem [M1, 12.2], at almost all points of the current boundary $\partial T$, T has a measure-theoretic exterior normal and the approximate tangent cone is a half-space of $\mathbf{R}^n\times\{0\}$. As in Almgren [Alm, V1.2(3)], gently pushing in that normal direction yields a smooth family $\Phi_t$ of diffeomorphisms for |t| small supported in a shrunken ball $B_1'' \subset B_1' \subset B_1$ such that $\Phi_0$ is the identity and the initial rate of change of the volume of T satisfies $dV/dt|_0 = 1$; let $A_0 = dA/dt|_0$ be the initial rate of change of boundary area. Choose $t_0 > 0$ so small that for all $-t_0 \leq t \leq t_0$, $.9 \leq dV/dt \leq 1.1$ and $|dA/dt| \leq |A_0| + 1$. View the $\Phi_t$ as a smooth family of diffeomorphisms of $\mathbf{R}^N$ supported in $B_1''\times\mathbf{R}^{N-n}$.

(2)     For $\alpha$ large, for all $-.9t_0 \leq \Delta V \leq .9t_0$, such diffeomorphisms of $T_\alpha$ alter volume by $\Delta V$ and increase the area of $\partial T$ at most $|\Delta V|(|A_0| + 1)/.9 = C_1|\Delta V|$.

In $M_\alpha$, let $B_\alpha = M_\alpha \cap (B_1'\times\mathbf{R}^{N-n})$ and let $\Phi_t^\alpha$ denote the G symmetrization of the $\Phi_t$, supported in $M_\alpha$ in $G(B_\alpha)$, which consists of the disjoint images of $B_\alpha$ under $G_1$ and distant copies thereof under other elements of G. Since the $T_\alpha$ are G-invariant, (2) still holds. We may assume that any $\delta_1/2$- ball in $M_\alpha$ is disjoint from the support of the $\Phi_t^\alpha$. Choose $\delta < \delta_1/2$ such that a $\delta$-ball in any $M_\alpha$ has volume less than $.9t_0$. We claim that for some $C_2 > 1/2\delta$ (independent of $\alpha$)

(3)     *in any $\delta$-ball B in $M_\alpha$, $\partial T_\alpha$ minimizes area plus $C_2|S-T_\alpha|$ in comparison with the restriction of other G-invariant surfaces $\partial S$ to B with the same boundary in $\partial B$.*

We may assume that S coincides with $\partial T_\alpha$ outside of G(B). Choose a ball $B_\alpha$ as above with $G(B_\alpha)$ disjoint from B and hence disjoint from G(B). By (2), obtain S' by altering S in $G(B_\alpha)$ so that it bounds net volume 0 with $\partial T_\alpha$ such that its area satisfies

$$|S'| \leq |S| + C_1|S-T_\alpha|_{G(B)} \leq |S| + |G|C_1|S-T_\alpha|_B.$$

Since $T_\alpha$ is isoperimetric among G-invariant regions,

$$|\partial T_\alpha| \leq |S'| \leq |S| + |G|\, C_1\, |S-T_\alpha|_B,$$

proving the claim (3). Next we claim that

(4)     *there is a positive lower bound (independent of $\alpha$) on the area of $\partial T_\alpha$ in any unit ball about one of its points in $M_\alpha$*

Consider a unit ball about a point p of $\partial T_\alpha$. For $r \leq \delta$, let $T_r = \partial T_\alpha \cap B(p, r)$ and let $g(r)$ denote the (n−1)D area of $T_r$. For almost all r, the boundary of $T_r$ has (n−2)D area at most $g'(r)$ (see [M1, Chapt. 9]). By an isoperimetric inequality for the manifold M [M1, Sect. 12.3], which is invariant under scaling of M in $\mathbf{R}^N$ and hence holds for each $M_\alpha$, this same boundary bounds a surface $T_r'$ of area

$$|T_r'| \leq \gamma g'(r)^{(n-1)/(n-2)},$$

for some isoperimetric constant γ; together with $T_r$, $T_r'$ bounds a region of mass or volume at most

$$\beta[g(r) + \gamma g'(r)^{(n-1)/(n-2)}]^{n/(n-1)} \leq C_3 g(r)^{n/(n-1)} + C_3 g'(r)^{n/(n-2)}$$

for some isoperimetric constant β and constant $C_3$. By claim (3),

$$|T_r| \leq |T_r'| + C_4 g(r)^{n/(n-1)} + C_4 g'(r)^{n/(n-2)} \leq \gamma g'(r)^{(n-1)/(n-2)} + C_4 g(r)^{n/(n-1)} + C_4 g'(r)^{n/(n-2)}.$$

Hence for some $C_5 > 0$,

$$g(r)\left(1 - C_5 g(r)^{1/(n-1)}\right) \leq C_5[g'(r)^{(n-1)/(n-2)} + g'(r)^{n/(n-2)}].$$

We may assume that the coefficient of $g(r)$ for any $r \leq 1/2C_2 \leq \delta$ is at least 1/2, since otherwise we have a lower bound as desired for $g(r)$. Let

$$\Omega = \{0 < r < \delta : g'(r) \leq 1\}.$$

We may assume that $|\Omega| \geq \delta/2$, since otherwise we immediately obtain a lower bounded as desired for $g(r)$. On Ω,

$$g(r) \leq C_6 \, g'(r)^{(n-1)/(n-2)}.$$

Since g is monotonically increasing, integration yields

$$g(r) \geq C_7 \, r^{n-1},$$

again yielding the desired lower bound for $g(1/2C_2)$, proving the claim (4).

By claim (4), there is a bound on the sums of the diameters of the components of $\partial R_\alpha$ after our scaling, so for the original $\partial R_\alpha$ before scaling, the sum of the diameters of the components of $\partial R_\alpha$ and hence of $R_\alpha$ goes to zero.

**3.3. Corollary.** *In a compact Euclidean orbifold, for small prescribed volume, an isoperimetric region is a small metric ball.*

*Proof.* By Proposition 3.2, a component of an isoperimetric region has small diameter. By Observation 2.3, it is a metric ball about the image of a fixed point. For such balls of volume V, perimeter is proportional to $V^{(n-1)/n}$. Since these functions are concave, a single ball is best.

*Remark.* Corollary 3.3 and proof also apply to orbifolds modeled on the sphere or on hyperbolic space.

**3.4. Theorem.** *In a compact orbifold O, for small prescribed volume, an isoperimetric region R is $C^0$ close to a small metric ball about a point where the group action has highest effective order. More specifically, if locally O is the quotient of a Riemannian geodesic ball M by a finite group G fixing its center, the preimage of R is $C^\infty$ close to a geodesic ball in M, and the order of G is the maximum over all such groups associated with O.*

*Remark.* In particular, in a smooth compact Riemannian manifold, for small volume, an isoperimetric region among regions invariant under a group G of isometries is the orbit of a nearly round ball.

*Proof.* The orbifold O is covered by quotients $M_i/G_i$ of Riemannian geodesic balls; we may assume that they still cover when shrunken by a factor of two. By Proposition 3.2, R has components of small diameter, each lying in the center of some $M_i/G_i$. Their preimages $R_i$ in $M_i$ are perimeter minimizing among $G_i$-invariant regions. The arguments of [MJ, Thm. 2.2] apply and show that each $R_i$ is $C^\infty$ close to a metric ball. Indeed, a limit of such $R_i$ with volume approaching 0 must be a minimizer in Euclidean space among $G_i'$-invariant regions, for some subgroup $G_i'$ of $G_i$, as in the proof of Proposition 3.2. By Observation 2.3, such a minimizer must be a round ball, and $R_i$ consists of an orbit under $G_i$ of $G_i'$-invariant balls. Easy estimates show that a single $G_i$-invariant ball is best. In particular, each component of R is $C^0$ close to a small metric ball (which is topologically a quotient of a ball).

      If there were more than one $R_i$, the second variation formula [MJ, 2.1] would imply instability, by uniformly expanding $R_1$ and shrinking $R_2$ (preserving enclosed volume). The final statement holds because for a small metric ball, the boundary area and volume satisfy

$$A \approx \beta V^{(n-1)/n}/|G_i|^{1/n},$$

where $A = \beta V^{(n-1)/n}$ in $\mathbf{R}^n$.

Isolated singularities in orbifolds are modeled on $\mathbf{R}^n$ modulo a finite subgroup of the orthogonal group with a single fixed point. The following consequence of [MR] applies to manifolds with more general isolated conical singularities.

**3.5. Proposition.** *Let $M^n$ be a smooth compact submanifold of $\mathbf{R}^N$ with a finite, positive number of conical singularities, where it is locally equivalent by a diffeomorphism of $\mathbf{R}^N$ to a cone with nonnegative Ricci curvature over a connected submanifold of $\mathbf{S}^{N-1}$. Then for small prescribed volume, an isoperimetric region is $C^\infty$ close to a metric ball about a singularity.*

*Proof.* Choose $0 < V_0 < 1$ by Bérard-Meyer [MR, Thm. 2.1] for the complement in M of unit balls about its singularities. A sequence of minimizers with volume $V_i \to 0$, scaled up by homothety about the $k^{th}$ singularity to volume $V_0$, converges to an isoperimetric region in the cone $C_k$ at that singularity. The region is nonempty for at most one such cone C. By [MR, Cor. 3.9], such an isoperimetric region in C must be a ball about the vertex. As in the proof of [MR, Thm. 2.2], no volume is lost to infinity. By [M1, Lemma 3.5], the (rescaled) sequence of minimizers has uniformly, weakly bounded mean curvature. By Allard's regularity theorem, the convergence to the ball is $C^{1,\alpha}$. By Schauder theory (see [M4, 3.3 and 3.5]), the convergence is $C^\infty$.

General singularities in orbifolds are modeled on $\mathbf{R}^n$ modulo a finite subgroup of the orthogonal group, where isoperimetric regions are metric balls (Cor. 3.3). The following consequence of [MR] extends Corollary 3.3.

**3.6. Proposition.** *Let $C^n$ be a cone with nonnegative Ricci curvature over a connected submanifold of the sphere in $\mathbf{R}^N$. Then isoperimetric regions in $C \times \mathbf{R}^m$ are metric balls about points in the product of the vertex with $\mathbf{R}^m$.*

*Proof.* By [MR, Cor. 3.9], isoperimetric regions in C are balls about the vertex. As noted in [MR, Cor. 3.12], by symmetrization horizontal slices in C of an isoperimetric region in $C \times \mathbf{R}^m$ are metric balls. Similarly, vertical slices in $\mathbf{R}^m$ are round balls. The minimizer in $C \times \mathbf{R}^m$ is determined by a generating curve in the first quadrant of the plane. It solves the same planar problem as the generator of the $\mathbf{O}(n) \times \mathbf{O}(m)$-invariant minimizer in $\mathbf{R}^{n+m}$, which is of course a sphere. (Namely, it minimizes $\int x^{n-1} y^{m-1} ds$ given $\int x^{n-1} y^{m-1} dA$.) Therefore the generating curve is a circle and the isoperimetric region is a metric ball.